%% file: purity_valuation_arXiv_v1.tex
\documentclass[10pt,a4paper, notitlepage]{article}
\usepackage{import} 
\import{}{preamble}

\bibliography{bibli.bib} 

\usepackage{libertinus}

	\usepackage{footnote}
	
	\newcommand{\inv}[1]{[\frac{1}{#1}]}
	\title{\textbf{\large{\uppercase{Flat Cohomological Purity for Syntomic Schemes over Valuation Rings}}}}
	\author{\uppercase{Arnab Kundu}}
\date{}
\usepackage{doi}
\begin{document}
	\maketitle
	\abstract{Grothendieck's cohomological purity predicts that the cohomology of a scheme is insensitive to removing a closed subscheme of sufficiently high codimension. In this article, we establish a form of flat cohomological purity over arbitrary (possibly infinite-rank) mixed-characteristic valuation rings $V$, thereby extending the theorem of Česnavičius–Scholze to the non-noetherian setting.\par More precisely, for a flat finite-type scheme over $V$ with local complete intersection fibres, we prove that the cohomology with coefficients in a commutative finite locally free group scheme remains unchanged after removing a closed subscheme satisfying a suitable fibrewise codimension condition; in particular, we obtain vanishing in low degrees and injectivity in the critical degree.\par As applications, we deduce purity results for local cohomology, for torsion in the Picard group, and for the Brauer group. In higher rank, our results yield sharper bounds than those previously obtained by Bhatt–Lurie and Madapusi–Mondal. The argument rests on recent advances in the structure theory of valuation rings.}
\setcounter{section}{-1}
\section{Flat cohomological purity and consequences}\label{scn:intro}
A basic principle of cohomological purity is that removing a closed subscheme of sufficiently large codimension does not affect cohomology. In the setting of flat cohomology with coefficients in commutative finite locally free group schemes, Česnavičius and Scholze \cite{ces-scholze} proved such a purity theorem for schemes with local complete intersection singularities.
\par The goal of this article is to extend this theorem to the non-noetherian setting of valuation rings. More precisely, we generalise \cite{ces-scholze} to obtain a version of flat cohomological purity for flat schemes of finite type with local complete intersection fibres over arbitrary mixed-characteristic valuation rings, including those of infinite rank, under a fibrewise codimension hypothesis. 
\par Schemes over (non-discrete) valuation rings have been extensively studied in the literature (see, for instance, \cite{gabber-ramero-foundations, fujiwara_kato_foundations_rigid_geometry}). Beyond their foundational role in $p$-adic geometry (e.g., rigid analytic and adic geometries; see \cite{huber96}), they provide a natural testing ground for extending cohomological statements beyond the noetherian setting. At the same time, their highly non-noetherian nature introduces substantial new difficulties.
\par This extension is therefore non-formal. Valuation rings are typically far from noetherian, often of infinite Krull dimension (see \ccite{datta_local_cohomology_valuation_ring}{\textsection7}), and classical proofs of cohomological purity relying on noetherian induction and dévissage arguments do not apply in this setting. Moreover, controlling codimension in families over the spectrum of such a ring requires a careful analysis of how the height of its primes interacts with the geometry of the fibres. As a result, new techniques are needed: our approach combines approximation by regular local rings, ramification-theoretic input for henselian valuation rings, and descent methods in the \textit{pro}-\textit{fppf} topology.
\par We fix some notation and conventions that will be applied throughout the article. We fix a prime number $p$. Unless mentioned otherwise, all cohomology groups are taken in the $\fppf$ topology. Given a closed subset $Z\subseteq X$ of a scheme $X$ and a sheaf $\scr{F}$ of abelian groups on $X$, the cohomology of $\scr{F}$ with supports in $Z$ is denoted by $H^q_Z(X,\scr{F})$. We recall that a \textit{syntomic morphism} is a flat morphism of finite presentation with locally complete intersection fibres (see \stacks{01ub}).
\par Our main result is the following, whose proof we outline at the end of this section.
\but{A}\label{thm:intro}
Let $\ca{O}$ be a valuation ring of possibly infinite rank that is flat over $\bb{Z}_{(p)}$, let $X$ be an $\ca{O}$-scheme that is the limit of a cofiltered system of syntomic $\ca{O}$-schemes with affine transition maps, and let $Z\hookrightarrow X$ be a closed subset whose associated topological space is noetherian (this condition is automatic if $\ca{O}$ is of finite rank). Suppose that $d\geqslant 1$ is an integer such that
\bn[(i)]
\item for a prime $\fr{p}\subset\ca{O}$ of finite height $h$, the fibre $Z_{\fr{p}}\subseteq X_{\fr{p}}$ is of codimension $\geqslant d-h$,
\item and, for all other primes $\fr{q}\subset\ca{O}$, the fibre $Z_{\fr{q}}\subseteq X_{\fr{q}}$ is of codimension $\geqslant d-1$.
\en
Then, given a commutative finite locally free $X$-group scheme $G$, we have \bud H^q_Z(X,G)=0,\text{ for each }q<d.\eud
\eut
The codimension condition in Theorem~\ref{thm:intro} reflects the interaction between the geometry of the fibres and the height of primes of $\ca{O}$: it interpolates between a codimension $\geqslant d$ condition on the generic fibre and weaker bounds over deeper points of $\spec(\ca{O})$.
\bur{1}[Finite étale case]
	For a commutative finite étale $X$-group scheme $G$, Theorem~\ref{thm:intro} reduces to a purity statement in étale cohomology. Indeed, by \cite[théorème~11.7]{brauerIII}, the flat and étale cohomology groups of such a group $G$ agree. In particular, in the case where $X$ is smooth, Theorem~\ref{thm:intro} extends results of Česnavičius \cite{ces-brauer} to the non-noetherian setting.
\eur
\vspace{0.2cm}\noindent\textbf{Relation to previous work:} As previously mentioned, Theorem~\ref{thm:intro} generalises the purity theorem of Česnavičius--Scholze \ccite{ces-scholze}{Theorem~7.1.2} (and Česnavičius \cite{ces-brauer}) from the noetherian setting to arbitrary mixed-characteristic valuation rings.
\par Purity results over valuation rings have been studied by Gabber--Ramero \ccite{gabber-ramero-foundations}{Proposition~11.4.7} and Guo--Liu \ccite{ning_liu_quasi-split}{Theorems~6.8 and 7.9} (similar results were also obtained in \ccite{phd-thesis}{Proposition~4.5}, see also \ccite{amar:unramified_case_gersten_conjecture}{Proposition~4.3}) for smooth schemes. In contrast, our results apply to schemes with local complete intersection fibres, and therefore extend these results to a broader class of singular schemes (see Remark~\ref{rem:smooth_case}).
\par Theorem~\ref{thm:intro} is also closely related to recent work of Bhatt–Lurie \cite{bhatt_lurie_absolute_prismatic_cohomology} and Madapusi–Mondal \cite{madapusi_mondal_perfect_f_gauges} on syntomic and flat cohomological purity (see Corollary~\ref{conj:bhatt_lurie}). However, already for rank $\geqslant 2$, our vanishing range improves on the known bounds in op.~cit.~(see Remark~\ref{rem:bhatt_lurie}).
\vspace{0.2cm}\newline\noindent\textbf{Applications.} Theorem~\ref{thm:intro} has several consequences for classical invariants. In particular, it yields purity results for local cohomology (Corollary~\ref{thm:Alocal}), for the torsion in the Picard group, and for the Brauer group (Corollary~\ref{conj:Gabber}).
\par In the case where $X=\spec(A)$ is local with closed point $z$ and $Z=\{z\}$, Theorem~\ref{thm:intro} immediately yields the following vanishing in local cohomology.
\buc{B}\label{thm:Alocal}
	Let $\ca{O}$ be a valuation ring of finite rank $r$ that is flat over $\bb{Z}_{(p)}$ and $A$ be a faithfully flat, local $\ca{O}$-algebra that is the localisation at a prime ideal of the colimit of a filtered system of syntomic $\ca{O}$-algebras. Suppose that the point corresponding to the maximal ideal $\fr{m}\subset A$ is $\ca{O}$-fibrewise of codimension $\geqslant d$ in $\spec(A)$. Then, given a commutative finite locally free $A$-group scheme $G$, we have	\bud H^q_{\fr{m}}(A,G)=0,\text{ for each }q<d+r.\eud 
\euc
\bur{2}[Sharpness of the bound]\label{rem:sharpness}
The computations in \ccite{ces-scholze}{Remark~7.2.6} show that already when $\ca{O}$ is a field (that is, when $r=0$), there are examples where the vanishing above does not hold for $q=d$. In particular, the bound is sharp.
\eur
The following, proved in \textsection\ref{sec:reduction_to_henselian_Case}, relies on the case of Theorem~\ref{thm:intro} where $G=\mu_n$, for some integer $n$. It establishes a relative version of Gabber's conjectures (\ccite{gabber_oberwolfach_report}{Conjectures~2 and 3}).
\buc{C}\label{conj:Gabber}
	Let $\ca{O}$ be a valuation ring of possibly infinite rank that is flat over $\bb{Z}_{(p)}$, let $X$ be the localisation of a syntomic $\ca{O}$-scheme, let $Z\hookrightarrow X$ be a closed subscheme of finite presentation. Suppose that $d\geqslant 3$ is an integer such that
	\bn[(i)]
	\item for a prime $\fr{p}\subset\ca{O}$ of finite height $h$, the fibre $Z_{\fr{p}}\subseteq X_{\fr{p}}$ is of codimension $\geqslant\max(1,d-h)$, and
	\item for all other primes $\fr{q}\subset\ca{O}$, the fibre $Z_{\fr{q}}\subseteq X_{\fr{q}}$ is of codimension $\geqslant d-1$.
	\en Then, we have \bd\label{vanishing:pic}\Pic(X)_{\tors}\cong\Pic(X\setminus Z)_{\tors}.\ed In particular, if $X$ is a local scheme and $Z$ is its closed point, then $\Pic(X\setminus Z)_{\tors}=0$. 
	\par Furthermore, if $d\geqslant 4$, then we even have that \bd\label{vanishing:br}H^2(X,\bb{G}_m)_{\tors}\twoheadrightarrow H^2(X\setminus Z,\bb{G}_m)_{\tors}.\ed
\euc
\bur{3}\label{rem:brauer_group} Let $\Br(X)$ be the \textit{Brauer group} of a scheme $X$ (called the \textit{Brauer}--\textit{Azumaya group} in \ccite{ct-skorobogatov_brauer_grothendieck_group_book}{Definition~3.1.3}). We recall that the equality $\Br(X)=H^2(X,\bb{G}_m)_{\tors}$ holds for any $X$ equipped with an ample line bundle (for instance, when $X$ is quasi-affine) thanks to the work of Gabber and de Jong \cite{de_jong_gabber_br} (see also \ccite{ct-skorobogatov_brauer_grothendieck_group_book}{\textsection4.2}). In particular, in this case, diagram \eqref{vanishing:br} can be rewritten as a surjection of Brauer groups \bd\tag{C.2$^{\prime}$}\Br(X)\twoheadrightarrow\Br(X\setminus Z).\ed
\eur
\par The following consequence of Theorem~\ref{thm:intro}, proved at the end of \textsection\ref{sec:reduction_to_henselian_Case}, recovers a special case of \ccite{bhatt_lurie_absolute_prismatic_cohomology}{Corollary~8.5.7} and \ccite{madapusi_mondal_perfect_f_gauges}{Theorem~8.3.6}. 
\par For a scheme $X$, let $R\Gamma_{\syn}(X,\bb{Z}_p(1))$ denote the syntomic complex as defined in \ccite{bhatt_lurie_absolute_prismatic_cohomology}{\textsection8.4}. For $r\geqslant 1$, we define the syntomic cohomology with torsion coefficients by $R\Gamma_{\syn}(X,\bb{Z}/p^r(1))\colonequals R\Gamma_{\syn}(X,\bb{Z}_p(1))\otimes_{\bb{Z}_p}^{\bb{L}}\bb{Z}/p^r\bb{Z}$ and we write its cohomology groups as $H^q_{\syn}(X,\bb{Z}/p^r(n))$. 
\buc{D}\label{conj:bhatt_lurie}
Let $\ca{O}$ be a valuation ring of rank one\footnote{Typical examples of a non-discrete valuation rings include the rings of integers in an algebraic closure of mixed-characteristic non-archimedean valued fields, and perfectoid valuation rings of rank one.} that is flat over $\bb{Z}_{(p)}$, let $X$ be an ind-syntomic $\ca{O}$-scheme, and let $Z\subset X$ be a closed subscheme. Suppose that $d\geqslant 1$ is an integer such that 
\hspace{2cm}\bun\setlength{\itemindent}{5.2cm}
\item $Z\inv{p}\hookrightarrow X\inv{p}$ is of codimension $\geqslant d$, and 
\item $Z_{p=0}\hookrightarrow X_{p=0}$ is of codimension $\geqslant d-1$.\footnote{We note that, given a scheme $S$ with an element $f\in\Gamma(S,\ca{O}_S)$, we shall denote the vanishing locus of $f$ by $S_{f=0}$ and non-vanishing locus by $S\inv{f}$.} 
\eun
Then, for any commutative finite locally free $X$-group scheme $G$, the restriction morphism \bd\label{claim:1:BL} H^q(X,G)\tir&H^q(X\setminus Z, G)\ed is an isomorphism for $q<d-1$ and an injection for $q< d$. 
In particular, if $d\geqslant 3$, then for any $r\geqslant 1$, we get an isomorphism of syntomic cohomology groups \bd\label{claim:2:BL} H^1_{\syn}(X,\bb{Z}/p^r(1))\iso H^1_{\syn}(X\setminus Z, \bb{Z}/p^r(1)).\ed
\euc
We emphasise that, unlike in \ccite{bhatt_lurie_absolute_prismatic_cohomology}{Corollary~8.5.7} and \ccite{madapusi_mondal_perfect_f_gauges}{Theorem~8.3.6}, it is not necessary in Corollary~\ref{conj:bhatt_lurie} that $Z$ lies inside $X_{p=0}$. 
\bur{4}[Higher rank valuation]\label{rem:bhatt_lurie}
Theorem~\ref{thm:intro} yields sharper bounds than \ccite{madapusi_mondal_perfect_f_gauges}{Theorem~8.3.6} when $\ca{O}$ is a valuation ring of finite rank $r\geqslant 2$. For example, consider  $X=\bb{A}^2_{\ca{O}}$ and $Z$ to be the $\ca{O}$-special fibre of $X$. Considering the maximal ideal $\fr{m}\subset\ca{O}$, we know that there exists an element $a\in\ca{O}$ such that $V(a)=\{\fr{m}\}\subseteq\spec(\ca{O})$ (this is a consequence of the prime avoidance lemma \stacks{00ds}, e.g., one may apply \ccite{phd-thesis}{Lemma~3.3(a)} to the unique prime $\fr{p}\subset\ca{O}$ of height $r-1$). In particular, even though $\spec(\ca{O})$ is of Krull dimension $2$, its closed point is cut out by a single equation $\{a=0\}$.
\par Set $X\inv{a}\colonequals X\setminus Z$. Then there is a long exact sequence \bud\textstyle 0\to\Gamma_Z(X,\ca{O}_X)\to\Gamma(X,\ca{O}_X)\xrightarrow{\mathrm{res}}\Gamma(X\inv{a},\ca{O}_{X\inv{a}})\to H^1_Z(X,\ca{O}_X)\to 0.\eud Since $\mathrm{res}$ is not surjective, it follows that $H^1_Z(X,\ca{O}_X)\neq 0$. Therefore, for a commutative finite locally free $X$-group scheme $G$, we \textit{cannot} deduce from \ccite{madapusi_mondal_perfect_f_gauges}{Theorem~8.3.6} that \bud H^1_Z(X,G)=0,\eud  even though this vanishing does follow from Theorem~\ref{thm:intro} for $d=2$.
\eur
\noindent\textbf{Overview of the proof of Theorem~\ref{thm:intro}.}
The proof proceeds in three main steps.
\bun
\item In \textsection\ref{sec:ind_regular_case}, we treat the case where $\mathcal{O}$ is ind-regular, reducing finally to the results of \cite{ces-scholze} using approximation by regular local rings and a limit argument.
\item In \textsection\ref{sec:henselian_Case}, in the henselian case, we pass to an absolute integral closure of $\mathcal{O}$ and use ramification-theoretic results of \cite{kato_thatte_upper_Ramification_valuation} together with descent in the \textit{pro}-\textit{fppf} topology to establish the desired vanishing.
\item Finally, we reduce the general case to the henselian case by an excision argument for flat cohomology, again relying on results of \cite{ces-scholze}.
\eun

\noindent\textbf{Motivation.} The expectation that such a purity statement should extend to valuation rings is motivated by the principle that such rings exhibit features analogous to regular local rings (see, e.g., \cite{kelly-morrow}), despite being typically highly non-noetherian. For instance, Zariski's uniformisation conjecture (see, for example, \ccite{phd-thesis}{Conjecture~2.1.1}) predicts that valuation rings should be ind-regular, and this is known in several important cases (see \cite{antieau_datta_derived_splinters, kelly-morrow, popescu_dimension_1_free_value_group_smooth, popescu_immediate_pure_transcendental_smooth}, also recalled here in Proposition~\ref{lem:pikachu-lemma}).
\par Moreover, valuation rings arising in mixed characteristic are well behaved from the perspective of $p$-adic geometry: by a theorem of Gabber–Ramero (see \ccite{gabber-ramero-almost}{Theorem~6.5.12}), any valuation ring that is flat over $\bb{Z}_p$ is $p$-quasi-syntomic (in the sense of \ccite{bhatt_lurie_absolute_prismatic_cohomology}{Definition~C.6}), and hence the same holds for the syntomic $\ca{O}$-scheme $X$. Such rings also arise naturally in approximation arguments in $p$-adic geometry (see, e.g., \cite{kerz-strunk-tamme, bhatt-mathew-arc}), and as local rings in several emerging Grothendieck topologies (see \cite{gabber_kelly_points, ehik}). These considerations suggest that cohomological purity phenomena should persist in this broader context, and Theorem~\ref{thm:intro} confirms this expectation.

\section{The ind-regular case}\label{sec:ind_regular_case}
In this section, our goal is to prove Proposition~\ref{lem:ind_regular_case}, which establishes the case of Theorem~\ref{thm:intro} where $\mathcal{O}$ either contains $\mathbb{Q}$ or its fraction field admits no nontrivial extensions of degree a power of the residue characteristic. In this situation, it is known that $\mathcal{O}$ can be expressed as a filtered colimit of regular local rings (see Proposition~\ref{lem:pikachu-lemma}). Consequently, in Proposition~\ref{lem:ind_regular_case}, we reduce Theorem~\ref{thm:intro} to the purity result of \ccite{ces-scholze}{Theorem~7.1.2} via a limit argument.
\par We begin with the following well-known lemma, which will be used throughout the remainder of the article.
\bl[\cite{sga4ii}]\label{lem:cohomology_commutes_with_limits}
	Let $X=\lim_{\alpha}X_{\alpha}$ be a filtered colimit in the category of quasi-compact and quasi-separated schemes with affine transition maps. Any commutative flat $X$-group scheme $G$ of finite presentation descends to $X_{\alpha}$ for $\alpha\gg 0$. Moreover, for any such $G$, the natural map is an equivalence \bud R\Gamma(X,G)\cong\colim_{\alpha}R\Gamma(X_{\alpha},G).\eud
\el
\bs
	This lemma is standard in the literature, and follows from the arguments in \ccite{sga4ii}{exposé~VI}. We include a sketch for convenience of the reader. 
	\par This can be proved in the same way as \stacks{09yq}. By \ccite{egaIV_3}{théorème~8.8.2(ii)} and \cite[Tags~\href{https://stacks.math.columbia.edu/tag/04AI}{04AI}, \href{https://stacks.math.columbia.edu/tag/0GS1}{0GS1} and \href{https://stacks.math.columbia.edu/tag/081E}{081E}]{stacks-project}, the commutative flat $X$-group scheme $G$ of finite presentation descends to a commutative flat $X_{\alpha}$-group scheme of finite presentation, for $\alpha\gg 0$. More generally, by a limit argument, the $\fppf$ site of $X$ can be expressed as a cofiltered limit of the \textit{fppf} sites of $X_{\alpha}$. Thus, we can conclude by following the arguments in \stacks{09yq}.
\es
Proposition~\ref{lem:pikachu-lemma} below appears in the work of Antieau--Datta \ccite{antieau_datta_derived_splinters}{Propositions~4.1.1 and~4.2.1, and Corollary~4.2.4} (see also Kelly--Morrow \ccite{kelly-morrow}{Proof of Theorem~3.1}; cf.~\cite{popescu_dimension_1_free_value_group_smooth, popescu_immediate_pure_transcendental_smooth}). Before proceeding, let us recall some necessary terminology. 
\par The \textit{residue characteristic exponent} $e$ of a local ring $\Lambda$ of residue characteristic $p$ is defined by $$e\colonequals\begin{cases} p, \text{ if }p>0\\ 1,\text{ otherwise}.\end{cases}$$ Given an integral domain $D$, we define the \textit{prime local subring} to be the smallest local ring contained as a subring $\Lambda\subseteq D$. 
\bun 
	\item[$\circ$] For example, if $A$ contains a field $\bb{F}$, then $\Lambda=\bb{F}$. 
	\item[$\circ$] Otherwise, $\Lambda=\bb{Z}_{(p)}$, where $p$ is the residue characteristic of $D$.
\eun
\par In the case where $\mathcal{O}$ contains a field of characteristic $0$, Proposition~\ref{lem:pikachu-lemma}, stated below, goes back to Zariski \cite{zariski_local_uniformisation}. In general, the statement verifies a special case of Zariski’s uniformisation conjecture (see, for example, \cite[Conjecture~2.1.1]{phd-thesis}).
\bp[\cite{antieau_datta_derived_splinters, kelly-morrow}]\label{lem:pikachu-lemma}
Let $\mathcal{O}$ be a valuation ring with prime local subring $\Lambda$ of residue characteristic exponent $p$. Assume that the fraction field $K$ of $\mathcal{O}$ admits no nontrivial finite extensions of degree a power of $p$. Then $\mathcal{O}$ can be written as a filtered union of regular local $\Lambda$-subalgebras of finite type.
\ep
\bp\label{lem:ind_regular_case}
	Let $\ca{O}$ be a valuation ring of finite rank that is a filtered colimit of regular local rings, \bd\label{diag:valuation_ring_lci}\ca{O}=\colim_{\delta\in\Delta}R_{\delta},\ed let $X$ be a syntomic $\ca{O}$-scheme, and let $Z\hookrightarrow X$ be a finitely presented closed subscheme. We suppose that $d\geqslant 1$ is an integer such that, for a prime $\fr{p}\subset\ca{O}$ of height $h$,
	\bud\text{the fibre $Z_{\fr{p}}\subseteq X_{\fr{p}}$ is of codimension $\geqslant d-h$.}\eud
	Then, given a commutative finite locally free $X$-group scheme $G$, we have \bud H^r_Z(X,G)=0,\text{ for each }r<d.\eud
\ep
\bs
	We consider the local-to-global spectral sequence from \ccite{sga4ii}{exposé~V, proposition~6.4} \bud H^p(X,\scr{H}^q_Z(X,G))\implies H^{p+q}_Z(X,G),\eud where $\scr{H}^q_Z(X,G)$ denotes the Zariski sheafification of presheaf $U\mapsto H^q_{Z\times_XU}(U,G)$. Therefore, to prove the claim, it reduces to show that $\scr{H}^{q}_Z(X,G)=0$ for all $q<d$. Taking the stalk at a point $z\in Z$, it suffices to demonstrate that \bd\label{diag:vansihing_local}\text{$H^q_{Z\cap\spec(A)}(A,G)=0$ for all $q<d$},\ed where $A$ is the local ring at $z\in Z$. However, since $\ca{O}$ is of finite rank, the topological space associated to the finite type $\ca{O}$-scheme $X$ is noetherian (see \stacks{0053}); as a consequence, the same holds for the topological space associated to $X'\colonequals\spec(A)$ (see \stacks{0052}). We argue by noetherian induction on the closed subset $Z'\colonequals Z\cap X'$, following the proof of \ccite{ces-scholze}{Theorem~7.1.2}, in order to reduce to the case when $Z'$ is a closed point. 
	\par By the long exact sequence of cohomology with supports, to show \eqref{diag:vansihing_local}, it suffices to verify that the morphism \bd\label{diag:vansihing_local_injection} H^q(X',G)\tir&H^q(X'\setminus Z', G)\ed is an injection for $q<d$ and an isomorphism for $q<d-1$. Letting $\eta\in Z'$ be a generic point, we consider an open neighbourhood $\eta\in U\subset X'$. Then, the Mayer--Vietoris sequence
	\bd\label{diag:Mayer--Vietoris}
		\cdots\tir&H^q(U\cup(X'\setminus Z'), G)\tir&H^q(U,G)\oplus H^q(X'\setminus Z', G)\tir& H^q(U\setminus Z',G)\tir&\cdots
	\ed
	proves that any class $[\xi]\in\ker(H^q(X',G)\to H^q(X'\setminus Z'), G)$ also vanishes over $U\setminus Z'$. However, by our hypothesis, the class $[\xi]$ restricted to the local ring $A_{\eta}$ of $U$ at $\eta$, which vanishes over  \bud\spec(A_{\eta})\setminus Z'=\spec(A_{\eta})\setminus\{\eta\},\eud must also vanish over $A_{\eta}$ itself. As a consequence, localising $U$ at $\eta$ and spreading out, we may assume that $[\xi]$ vanishes over $U$ already. Hence, this $[\xi]$ vanishes in the complement of $Z_U\colonequals Z\setminus U\subseteq X'$. Consequently, if $q<d$, by the induction hypothesis, this class must vanish on $X'$. This demonstrates that \eqref{diag:vansihing_local_injection} is an injection for $q<d$. In a similar vein, for any $q<d-1$, any class $[\xi]\in H^q(X'\setminus Z', G)$ extends to a neighbourhood $\eta\in U$. Therefore, by the induction hypothesis, this $[\xi]$ extends all over to $X'$. This concludes the induction. As a result, it suffices to demonstrate below the vanishing \eqref{diag:vansihing_local} in the special case when $Z=\{z\}$ is the closed point.
	\par We suppose that $z$ lies over a point $\fr{p}$ of height $r$ in $\spec(\ca{O})$, localising at which we may assume that $\spec(A)\to\spec(\ca{O})$ is surjective. Furthermore, by localising each $R_{\delta}$ at the image $\fr{p}_{\delta}$ of $\fr{p}$ under the morphism $\spec(\ca{O})\to\spec(R_{\delta})$, we may preserve the colimit \eqref{diag:valuation_ring_lci}, which induces an isomorphism of schemes \bud\textstyle\spec(\ca{O})\cong\lim_{\delta\in\Delta}\spec(R_{\delta}).\eud Consequently, as $\spec(\ca{O})$ is a finite topological space, after possibly refining $\Delta$, we may assume that $\spec(\ca{O})$ embeds inside $\spec(R_{\delta})$ for each $\delta\in\Delta$. In particular, each $R_{\delta}$ is of Krull dimension $\geqslant r$, or equivalently, the maximal ideal of each $R_{\delta}$ is of height $\geqslant r$. 
	\par Since $\ca{O}$ is coherent (see \stacks{0ewv}), the flat $\ca{O}$-scheme $X$ of finite type is of finite presentation. In particular, by \ccite{egaIV_3}{théorème~8.8.2(ii)} and \cite[Tags~\href{https://stacks.math.columbia.edu/tag/0C3L}{0C3L} and \href{https://stacks.math.columbia.edu/tag/07RR}{07RR}]{stacks-project}, possibly by modifying $\Delta$, we may assume that $X\to\spec(\ca{O})$ descends to a syntomic surjective morphism $X_{\delta}\to\spec(R_{\delta})$ for each $\delta\in\Delta$. 
	Furthermore, by \ccite{egaIV_3}{théorème~8.8.2(ii)} and \cite[Tags~\href{https://stacks.math.columbia.edu/tag/04AI}{04AI} and \href{https://stacks.math.columbia.edu/tag/06AC}{06AC}]{stacks-project}, and possibly modifying $\Delta$ again, we may assume that $G$ descends to a finite locally free $X_{\delta}$-scheme, also denoted by $G$ for simplicity, for each $\delta\in\Delta$. Finally, by \cite[Tags~\href{https://stacks.math.columbia.edu/tag/0GS1}{0GS1} and \href{https://stacks.math.columbia.edu/tag/081E}{081E}]{stacks-project}, it follows that this $X_{\delta}$-scheme $G$ is a commutative group scheme for each $\delta\in\Delta$.
	\par We define $A_{\delta}$ to be the local ring of $X_{\delta}$ at the image $z_{\delta}$ of $z$ along the morphism $X\to X_{\delta}$. By construction, we have that \bd\label{diag:colimit} A\cong\colim_{\delta\in\Delta} A_{\delta},\ed where each transition map is a local morphism. By a limit argument and a further modification of $\Delta$, since $z$ is of codimension $\geqslant d-r$ in the $\ca{O}$-fibre of $\spec(A)$ over $\fr{p}$, it follows that $z_{\delta}$ is also of codimension $\geqslant d-r$ in the $R_{\delta}$-fibre of $\spec(A_{\delta})$ over $\fr{p}_{\delta}$. 
	\par However, we note that each $A_{\delta}$ is the local ring of a syntomic scheme over the regular local ring $R_{\delta}$, and is therefore a noetherian local ring. Moreover, by the local structure of syntomic morphisms (see \stacks{01ue}), each $A_{\delta}$ is a complete intersection ring. Since $R_{\delta}$ is of Krull dimension $\geqslant r$ and $z_{\delta}$ is of codimension $\geqslant d-r$ in its fibre, by, for example, the catenary property of $A_{\delta}$ and the $R_{\delta}$-flatness of $A_{\delta}$ (\stacks{00hs}), the Krull dimension of $A_{\delta}$ is $\geqslant d$. In particular, thanks to \ccite{ces-scholze}{Theorem~7.1.2}, we obtain \bud H^q_{\{z_{\delta}\}}(A_{\delta},G)=0\text{ for }q<d.\eud 
	As a result, the definition of cohomology with supports along with \eqref{diag:colimit} and Lemma~\ref{lem:cohomology_commutes_with_limits}, we have that \bud H^q_{\{z\}}(A,G)=0\text{ for }q<d,\eud as required. This terminates the proof.  
\es
\section{The henselian case: reduction to absolute integral closure}\label{sec:henselian_Case}
In this section, our goal is to prove Proposition~\ref{lem:reduction_to_ind_regular}, which establishes the case of Theorem~\ref{thm:intro} where $\ca{O}$ is henselian. We reduce the proof to the case treated by Proposition~\ref{lem:ind_regular_case}. To this end, we follow the arguments in \ccite{ces-scholze}{Lemma~4.1.12}, and, using flat descent, we pass to a suitable integral valuation ring extension $\ca{O}'$ of $\ca{O}$. 
\par This passage is delicate: beyond the discrete case, even when the induced extension of fraction fields is finite, this $\ca{O}'$ is rarely of finite type over $\ca{O}$, as shown by their ramification theory (see Proposition~\ref{prop:kato_thatte}). As a consequence, this reduction requires the use of \textit{pro}-\textit{fppf} descent (recalled in Definition~\ref{defn:prof-fppf}).
\par We begin by recalling the following definition from \ccite{ces-scholze}{\textsection4.1.5}, which is required to state Proposition~\ref{prop:kato_thatte}.
\bdf\label{defn:prof-fppf} 
	A faithfully flat morphism $Y\to X$ of quasi-compact and quasi-separated schemes is called a \textit{pro-fppf cover} if for each $y\in Y$, there exists an affine open $y\in\spec(B)\subseteq Y$ whose image lies in an affine open $\spec(A)\subseteq X$ such that the induced morphism $A \to B$ is a filtered colimit of flat $A$-algebras of finite presentation.
\edf
We will need the following properties of \textit{pro}-\textit{fppf} topology in the proof of Proposition~\ref{prop:kato_thatte}.
\br\label{rem:pro-fppf}
	\bn[(i)]
		\item\label{pt:i:rem:pro-fppf} By \ccite{ces-scholze}{\textsection2.3.3}, $\it{pro}\text{-}\it{fppf}$ covers of quasi-compact and quasi-separated schemes are stable under base change and composition. Any such \textit{pro}-\textit{fppf} cover is, in particular, an \textit{fpqc} cover.
		\item Given a commutative finite locally free group scheme $G$ over a quasi-compact and quasi-separated scheme $X$, the cohomology functor $R\Gamma(-,G)$ satisfies descent in the \textit{pro}-\textit{fppf} topology over $X$. Indeed, thanks to \ccite{ces-scholze}{Theorem~5.5.2}, this $R\Gamma(-,G)$ satisfies hyperdescent along any faithfully flat extension of affine schemes $\spec(A')\to\spec(A)$ defined over $X$. Consequently, it satisfies descent along any faithfully flat cover of quasi-compact and quasi-separated $X$-schemes. Therefore, by \eqref{pt:i:rem:pro-fppf}, this $R\Gamma(-,G)$ is a sheaf for the \textit{pro}-\textit{fppf} topology on $X$, and by taking fibres, the same holds for the cohomology $R\Gamma_Z(-,G)$ with supports in a closed subscheme $Z\subseteq X$.  
	\en
\er
\par The following result is essentially \ccite{kato_thatte_upper_Ramification_valuation}{Theorem~6.0.2} (cf.~\ccite{kerz-strunk-tamme}{Proposition~A.10} and \cite{popescu_neron, popescu_algebraic_complete_intersection, popescu_char_p_complete_intersection, popescu_free_value_group_complete_intersection, popescu_mixed_char_torsion_free_value_group_complete_intersection}).
\bp[\cite{kato_thatte_upper_Ramification_valuation}]\label{prop:kato_thatte}
The integral closure $\tilde{\ca{O}}$ of a henselian valuation ring $\ca{O}$ with fraction field $K$ in a (possibly infinite) Galois extension $\tilde{K}/K$ is a faithfully flat valuation ring extension of $\ca{O}$, that can be written as a filtered colimit \bud\tilde{\ca{O}}=\colim R,\eud of syntomic $\ca{O}$-subalgebras $R\subseteq\tilde{\ca{O}}$ that are free of finite rank as $\ca{O}$-modules. In particular, the induced morphism $\spec(\ca{\tilde{O}})\to\spec(\ca{O})$ is a pro-fppf cover and induces a homeomorphism on underlying topological spaces.
\ep
\bs
As $\mathcal{O}$ is a valuation ring, the algebra $\tilde{\mathcal{O}}$ is $\mathcal{O}$-flat because it is torsion-free as an $\ca{O}$-module. By the going-up property for integral morphisms \stacks{00gu}, the map $\mathcal{O}\hookrightarrow\tilde{\mathcal{O}}$ is faithfully flat. Moreover, each fibre over $\mathfrak{p}\subset\mathcal{O}$ is a $0$-dimensional integral domain, hence a field (\cite[Tags~\href{https://stacks.math.columbia.edu/tag/02JK}{02JK} and \href{https://stacks.math.columbia.edu/tag/00GS}{00GS}]{stacks-project}). Thus $\spec(\tilde{\mathcal{O}})\to\spec(\mathcal{O})$ is bijective, and therefore a homeomorphism because both spectra are totally ordered.
	\par In order to show the rest of the claims, it suffices to assume that the extension $\tilde{K}/K$ is moreover finite. Indeed, any Galois extension is a filtered union of finite Galois extensions, and filtered colimits commute. Consequently, the required claims are the content of \ccite{kato_thatte_upper_Ramification_valuation}{Theorem~6.0.2}.
\es
\br[Purely inseparable case]\label{rem:purely_insep}
	The conclusion of Proposition~\ref{prop:kato_thatte} is known in several cases where the extension $\tilde{K}/K$ is not Galois, including when it is purely inseparable. When $\ca{O}$ has rank $1$, the latter case follows from arguments of Gabber (see \ccite{kerz-strunk-tamme}{\textsection A.2}). In general, one expects this to follow from Popescu’s work \cite{popescu_algebraic_complete_intersection}.
\er
Let $X$ be a scheme. Given a family $\scr{U}$ of morphisms with target $X$ and a presheaf $\scr{F}$ of abelian groups over $X$, let \bd\label{diag:cech_complex}\check{\ca{C}}^0(\scr{U},\scr{F})\to\check{\ca{C}}^1(\scr{U},\scr{F})\to\check{\ca{C}}^2(\scr{U},\scr{F})\to\check{\ca{C}}^3(\scr{U},\scr{F})\to\cdots\ed denote the \textit{\v{C}ech complex} of $\scr{F}$ with respect to $\scr{U}$ (the individual terms are defined in \stacks{03ol}).
\bp\label{lem:reduction_to_ind_regular}
	Let $\ca{O}$ be a henselian valuation ring of finite rank that is flat over $\bb{Z}_{(p)}$, let $X$ be a syntomic $\ca{O}$-scheme, and let $Z\hookrightarrow X$ be a closed subscheme. Suppose that $d\geqslant 1$ is an integer such that
\bn[(i)]
\item\label{eqref:i} for a prime $\fr{p}\subset\ca{O}$ of finite height $h$, the fibre $Z_{\fr{p}}\subseteq X_{\fr{p}}$ is of codimension $\geqslant d-h$, and
\item\label{eqref:ii} for all other primes $\fr{q}\subset\ca{O}$, the fibre $Z_{\fr{q}}\subseteq X_{\fr{q}}$ is of codimension $\geqslant d-1$.
\en
Then, given a commutative finite locally free $X$-group scheme $G$, we have \bd\label{prop:finite_extension_case} H^r_Z(X,G)=0,\text{ for each }r<d.\ed
\ep
\bs
	Let $\tilde{K}$ be a Galois extension of the fraction field $K$ of $\ca{O}$ and let $\tilde{\ca{O}}\subset\tilde{K}$ be the integral closure of $\ca{O}$. Given that $\ca{O}$ is henselian, by Proposition~\ref{prop:kato_thatte}, it follows that $\ca{O}\hookrightarrow\tilde{\ca{O}}$ is a faithfully flat extension of valuation rings. We suppose that $p$ is the residue characteristic of $\ca{O}$, and henceforth, assume that $\tilde{K}$ is the algebraic closure of $K$. Consider the base changes $\tilde{Z}\colonequals Z\times_{\ca{O}}\spec(\tilde{\ca{O}})\hookrightarrow \tilde{X}\colonequals X\times_{\ca{O}}\spec(\tilde{\ca{O}})$, inducing a cartesian diagram of schemes \bd\label{diag:base_Change}\tilde{Z}\arrow[r,hook]\arrow[d]&\tilde{X}\arrow[d]\\Z\arrow[r,hook]&X.\ed Since flatness as well as the property of being syntomic is preserved by any base change and flat morphism of rings satisfy going down (\stacks{00hs}), we verify that $\tilde{Z}\hookrightarrow\tilde{X}$ satisfies the fibre codimension property over any prime $\tilde{\fr{p}}\subset\tilde{\ca{O}}$ for the same $d$. Indeed, this can be checked fibrewise using the homeomorphism $\spec(\tilde{\ca{O}})\to\spec(\ca{O})$ (see Proposition~\ref{prop:kato_thatte}). Thus, by Propositions~\ref{lem:pikachu-lemma} and \ref{lem:ind_regular_case}, we get \bud H^r_{\tilde{Z}}(\tilde{X},G)=0,\text{ for each }r<d.\eud In particular, to prove \eqref{prop:finite_extension_case}, it suffices to show that \bd\label{diag:injection} H^r_{Z}(X,G)\hookrightarrow H^r_{\tilde{Z}}(\tilde{X},G),\text{ for each }r<d.\ed
	We establish this by following the arguments in \ccite{ces-scholze}{Lemma~4.1.12}. Again, by Proposition~\ref{prop:kato_thatte}, the morphism $\spec(\tilde{\ca{O}})\to\spec(\ca{O})$ is a \textit{pro}-\textit{fppf} cover; hence, by Remark~\ref{rem:pro-fppf}, the same follows for the morphism $\tilde{X}\to X$. 
	\par Before proceeding, let us introduce some notation for the rest of the proof. For any $X$-scheme $T$ and any integer $q$, we write $H^q_Z(T,-)$ for $H^q_{Z\times_XT}(T,-)$. Furthermore, for any integer $m\geqslant 1$, the $m$-fold product $T\times_X\cdots\times_XT$ will be denoted by $T^m$.
	\par Since \textit{fppf} cohomology (with supports) has \textit{pro}-\textit{fppf} descent, we consider the Čech-to-derived spectral sequence \stacks{03ow} (also called the Cartan--Leray spectral sequence in \ccite{sga4ii}{exposé~V, corollaire~3.3}) with respect to the \textit{pro}-\textit{fppf} cover $\scr{U}\colonequals\{\tilde{X}\to X\}$ \bud E_1^{p,q}\colonequals\check{\ca{C}}^p(\scr{U}, H^q_{Z}(-,G))\Rightarrow H^{p+q}_Z(X,G),\eud where $H^q_Z(-,G)$ denotes the presheaf $U\mapsto H^q_{Z}(U,G)$ for each $X$-scheme $U$ and for each $p$ (assuming the notation of \eqref{diag:cech_complex}). 
	Thus, to prove \eqref{diag:injection} for $r=n$, it suffices to verify that $E_1^{p,n-1}=0$ for each $p\geqslant 0$. Indeed, in this case, we obtain that $H^n_Z(X,G)\cong E_{\infty}^{0,n}\hookrightarrow E_1^{0,n}=H^n_Z(\tilde{X},G)$. Thus, it remains to show that for each $m\geqslant 1$,  \bd\label{diag:product_vanishing}\displaystyle H^{n-1}_{Z}(\tilde{X}^m,G)=0.\ed Since the property of being integral as well as syntomic are preserved under base changes in addition to compositions, each product $\tilde{X}^m$ is flat and integral over $X$ as well as syntomic over $\tilde{\ca{O}}$. Consequently, for each $p\geqslant 0$, the $\tilde{\ca{O}}$-closed immersion $\tilde{Z}^m\hookrightarrow\tilde{X}^m$ produces a cartesian diagram of schemes analogous to \eqref{diag:base_Change}. In particular, as before, checking fibrewise over each prime $\tilde{\fr{p}}\subset\tilde{\ca{O}}$, we can verify that $\tilde{Z}^m\hookrightarrow\tilde{X}^m$ satisfies the required codimension hypothesis of this proposition for the same $d$ over $\tilde{\ca{O}}$. Consequently, Proposition~\ref{lem:ind_regular_case} yields the necessary vanishing \eqref{diag:product_vanishing}.
\es
\br[Mixed-characteristic assumption]
\bun
\item[$\circ$] We note that the proof of Proposition~\ref{lem:reduction_to_ind_regular} only requires that $\tilde{K}$ admit no nontrivial extensions of degree a power of $p$. In particular, one could instead take $\tilde{K}$ to be the extension generated by all finite extensions of $K$ of degree a power of $p$.
\item[$\circ$] In mixed characteristic, any such extension $\tilde{K}/K$ is automatically separable, so Proposition~\ref{prop:kato_thatte} applies. In particular, if the conclusion of Proposition~\ref{prop:kato_thatte} were known for purely inseparable extensions (see Remark~\ref{rem:purely_insep}), then Theorem~\ref{thm:intro} would extend to the equal-characteristic $p$ case. We plan to return to this question in future work.
\eun
\er
\section{Reduction to the henselian case: excision}\label{sec:reduction_to_henselian_Case}
In this section, we finish the proof of Theorem~\ref{thm:intro}. We first reduce the statement to the case where $\mathcal{O}$ is of finite rank (Theorem~\ref{thm:fd_case}), and then further to the henselian case treated in Proposition~\ref{lem:reduction_to_ind_regular}. The key input is the following excision result for \textit{fppf} cohomology, recalled from \ccite{ces-scholze}{Theorem~5.4.4} below.
\bp[\cite{ces-scholze}]\label{prop:excision}
Given a flat extension of rings $A\to A'$ , a finitely generated ideal $I\subseteq A$ satisfying that $A/I\iso A'/IA'$ and a commutative finite locally free $A$-group scheme $G$, we have 
\bud R\Gamma_{I}(A,G)\diso & R\Gamma_{IA'}(A', G).\eud
\ep
\bs
	The flatness of $A\to A'$ ensures that the derived tensor product $A/I\otimes^{\bb{L}}A'$ coincides with the usual tensor product $A'/IA'$. The claim then is the content of \ccite{ces-scholze}{Theorem~5.4.4}.
\es
\bt\label{thm:fd_case}
Let $\ca{O}$ be a valuation ring of finite rank that is flat over $\bb{Z}_{(p)}$, let $X$ be the localisation of a syntomic $\ca{O}$-scheme and let $Z\hookrightarrow X$ be a closed subset. Suppose that $d\geqslant 1$ is an integer such that
\bud\text{for a prime $\fr{p}\subset\ca{O}$ of height $h$, the fibre $Z_{\fr{p}}\subseteq X_{\fr{p}}$ is of codimension $\geqslant d-h$.}\eud
Then, given a commutative finite locally free $X$-group scheme $G$, we have \bud H^q_Z(X,G)=0,\text{ for each }q<d.\eud
\et
\bs
Since the spectrum of $\ca{O}$ is a finite topological space, the underlying topological space of the finite type $\ca{O}$-scheme $X$ is noetherian. As a consequence, if $\mathscr{I}\subseteq\mathcal{O}_X$ is a quasi-coherent sheaf of ideals vanishing on $Z\subseteq X$, then writing $\mathscr{I}$ as a filtered colimit of its finitely generated subideals and using the noetherian property, we conclude that $Z$ is defined by a finitely generated ideal. In particular, $Z$ admits a finitely presented closed subscheme structure.
\par Arguing as in the proof of Proposition~\ref{lem:ind_regular_case}, the local-to-global spectral sequence \ccite{sga4ii}{exposé~V, proposition~6.4} reduces the problem to showing that \bd\label{diag:vanishing:redn_to_henselian}\text{$H^q_{Z\cap\spec(A)}(A,G)=0$ for all $q<d$},\ed where $A$ is the local ring at some point $z\in Z$. Consequently, without loss of generality, we may assume that $X=\spec(A)$ and $Z=Z\cap\spec(A)$. Note that this does not change the fact that $X$ is noetherian because it is the subspace of a noetherian topological space.
\par In a similar vein as the proof of Proposition~\ref{lem:ind_regular_case}, choose an open $U\subset X$ containing a generic point $\eta\in Z$. Localising $U$ at $\eta$ and then spreading out, a noetherian induction argument using the Mayer--Vietoris sequence \eqref{diag:Mayer--Vietoris} reduces the proof to the case of \eqref{diag:vanishing:redn_to_henselian} where $Z=\{z\}$ is the closed point.
\par Suppose that $z$ lies over a point of height $r$ in $\spec(\mathcal{O})$. After localising at this point, we may assume that $X\to\spec(\mathcal{O})$ is surjective and that $\mathcal{O}$ is of rank $r$. By Proposition~\ref{prop:excision}, the vanishing \eqref{diag:vanishing:redn_to_henselian} in the case $Z=\{z\}$ is equivalent to \bd\label{diag:vansihing_henselian}\text{$H^q_{\{z^h\}}(A^h,G)=0$ for all $q<d$,}\ed where $A^h$ is the henselisation of $A$ with closed point $z^h\in X^h\colonequals\spec(A^h)$. Let $\ca{O}^h$ be the henselisation of the valuation ring $\ca{O}$ (see \stacks{0ask}). By prime avoidance lemma \stacks{00ds}, there exists an element $a\in\ca{O}^h$ in the maximal ideal $m^h\subset\ca{O}^h$ so that $V(a)=\{m^h\}$. Arguing as above, the closed point $z\in X$ is defined by a finitely generated ideal $I\subset A$. Replacing $I$ by $I+aA$ if necessary, we may assume that $a\in I$. In particular, since $A^h/\fr{m}A^h\cong A/\fr{m}$, the closed point $z^h\in X^h$ is defined by the ideal $IA^h$.
\par We claim that vanishing \eqref{diag:vansihing_henselian} will follow, if we establish that $X^h$ is an ind-syntomic $\ca{O}^h$-scheme such that $z^h$ is of codimension $\geqslant d-r$ in the $\ca{O}^h$-special fibre of $X^h$. 
Indeed, in this case, we can write $A$ as a filtered colimit of syntomic $\ca{O}^h$-algebras \bud A^h=\colim_{\delta\in\Delta}A_{\delta}.\eud Since $IA^h$ is finitely generated, we may assume that there is an ideal $I_{\delta}\subset A_{\delta}$ such that $I_{\delta}A^h=IA^h$ for $\delta\gg 0$. Possibly modifying $\Delta$, for each $\delta\in\Delta$, we may therefore define \bud Z_{\delta}\colonequals V(I_{\delta})\hookrightarrow X_{\delta}\colonequals\spec(A_{\delta}).\eud Since $a\in I_{\delta}$, it follows that each $Z_{\delta}$ lies in the $\ca{O}^h$-special fibre $(X_{\delta})_{\kappa}$ of $X_{\delta}$. Whereas, since ${z^h}=V(IA^h)$ is of codimension $\geqslant d-r$, for $\delta\gg 0$, we may assume that $Z_{\delta}=V(I_{\delta})$ is of codimension $\geqslant d-r$ in $(X_{\delta})_{\kappa}$. Indeed, since the special $\ca{O}^h$-fibre $X^h_{\kappa}$ of $X^h$ satisfies that $X^h_{\kappa}\cong\lim_{\delta\in\Delta} (X_{\delta})_{\kappa}$, this can be proved by picking and then descending to each $(X_{\delta})_{\kappa}$ a finite chain $C$ of length $d-r$ of pairwise strict specialisations of points starting at a generic point of $X^h_{\kappa}$ and ending at $z^h$. Moreover, modifying $\Delta$ further, we may assume that $G$ is actually defined over $A_{\delta}$ for each $\delta\in\Delta$. Therefore, thanks to Proposition~\ref{lem:reduction_to_ind_regular}, we get\bud H^r_{Z_{\delta}}(X_{\delta},G)=0,\text{ for each }r<d.\eud Therefore, by a limit argument using Lemma~\ref{lem:cohomology_commutes_with_limits}, we get the required vanishing \eqref{diag:vansihing_henselian}.

\par Consequently, it suffices to show that $X^h$ is an ind-syntomic scheme over $\ca{O}^h$ such that $z^h$ is of codimension $\geqslant d-r$ in $X^h_{\kappa}$. First, we note that the following morphism of topological spaces \bud h\colon\spec(\ca{O}^h)\to\spec(\ca{O})\eud is a homeomorphism (follows, for example, from \stacks{06lk} or the fact that $\ca{O}^h$ is ind-étale and hence of relative dimension $0$ over the valuation ring $\ca{O}^h$). We can therefore verify that $z^h$ is of codimension $\geqslant d-r$ in $X^h_{\kappa}$. Thus, it remains to show that $X^h$ is an ind-syntomic scheme over $\ca{O}^h$.
\begin{tpic}
	\node (Ah) at (0,1.6) {$X^h$};
	\node (A) at (3,0.8) {$X$};
	\node (Oh) at (0,0) {$\spec(\ca{O}^h)$};
	\node (O) at (3,-0.8) {$\spec(\ca{O})$};
	\path[->,font=\scriptsize,>=angle 90] (Ah) edge node[above] {$h_A$} (A) edge node[left] {$g$} (Oh) (A) edge node[right] {$f$} (O) (Oh) edge node[below] {$h$} (O);
\end{tpic}To see this, note that since $h$ is bijective, the morphism $g$ defined above is also surjective. In particular, the $\mathcal{O}^h$-algebra $A^h$ is torsion-free as an $\mathcal{O}^h$-module, and hence flat. Moreover, the composite $h\circ g = f\circ h_A$ is ind-syntomic, being the composition of the ind-étale morphism $h_A$ with the syntomic morphism $f$ (see \stacks{01ub}). Since $h$ is bijective, the fibres of $g$ identify with those of $h\circ g$. It then follows from \stacks{01uf} that $f_A$ is ind-syntomic, as claimed.
\es
\bs[Proof of Theorem~\ref{thm:intro}]
	Arguing as in the proof of Theorem~\ref{thm:fd_case}, since $Z$ is noetherian as a topological space, we first reduce to the case where $X=\spec(A)$ is the local ring at a point $z$ of an $\mathcal{O}$-algebra that is a colimit of filtered system of syntomic $\mathcal{O}$-algebras, with $Z=\{z\}$. By a limit argument, therefore, we may moreover assume that $A$ is actually the local ring of a syntomic $\mathcal{O}$-algebra. If $\mathcal{O}$ has finite rank, the claim follows from Theorem~\ref{thm:fd_case}. Thus, it remains to treat the infinite rank case.
	\par Given an arbitrary valuation ring $\ca{O}$, we know that it can be expressed as a filtered colimit \bud\ca{O}=\colim_{\delta\in\Delta}\ca{O}_{\delta}\eud of valuation rings of finite rank (this is well known in the literature, see, for instance, \ccite{bhatt-mathew-arc}{Lemma~2.22}).
	Concretely, one may write the fraction field of $\ca{O}$ as a filtered union of subfields $K_\delta$ of finite transcendence degree over the prime subfield, and then set $\ca{O}_{\delta} \coloneqq \mathcal{O} \cap K_{\delta}$, which are valuation rings of finite rank. 
	\par Since $X$ is flat and essentially of finite type over the coherent ring $\mathcal{O}$ (see \stacks{0ewv}), it is in fact essentially of finite presentation. Consequently, after possibly refining $\Delta$, we may assume that $X$ descends to $X_{\delta}=\spec(A_{\delta})$, where $A_{\delta}$ is the local ring of a syntomic $\mathcal{O}_{\delta}$-scheme at a point $z_{\delta}$ (see \stacks{01zf}), with $z_{\delta}$ the image of $z$. This yields a cartesian square \bd\label{diag:cartesian_square} \{z\}\arrow[r,hook]\arrow[d]&X\arrow[d]\\ \{z_{\delta}\}\arrow[r,hook]& X_{\delta}.\ed
	\par We claim that it is enough to show that top morphism of \eqref{diag:cartesian_square} satisfies the hypothesis of Theorem~\ref{thm:fd_case} for the same $d\geqslant 1$. Indeed, in that case, we obtain that \bud H^r_{\{z_{\delta}\}}(X_{\delta},G)=0,\text{ for each }r<d,\eud and the desired vanishing follows by passing to the limit via Lemma~\ref{lem:cohomology_commutes_with_limits}.
	\par It therefore remains to check the fibrewise codimension condition. Let $\mathfrak{p}\subset\mathcal{O}$ (resp., $\mathfrak{p}_{\delta}\subset\mathcal{O}_{\delta}$) be the prime lying below $z$ (resp., $z_{\delta}$), of height $h$ (resp., $h_{\delta}$). After localising $\mathcal{O}$ at $\mathfrak{p}$, we may assume that $X\to\spec(\mathcal{O})$ is surjective. By \stacks{07rr}, it follows that $X_{\delta}\to\spec(\mathcal{O}_{\delta})$ is surjective for $\delta\gg 0$, and we may assume this holds for all $\delta$. Suppose further that $h_{\delta}\geqslant 1$ for $\delta\gg 0$, because otherwise $h=0$.
	\par Since $\mathcal{O}$ has infinite rank, by the surjectivity assumption on $X\to\spec(\ca{O})$, we have $h=\infty$. By hypothesis, the fibre of the top morphism in \eqref{diag:cartesian_square} over $\mathfrak{p}$ has codimension $\geqslant d-1$. Since $X\cong\lim_{\delta\in\Delta}X_{\delta}$, the same holds for the fibre of $X_{\delta}$ over $\mathfrak{p}_{\delta}$ for $\delta\gg 0$. Indeed, one may descend to $X_{\delta}$ a finite chain $C$ of length $d-1$ of pairwise strict specialisations of points in $X$ from a generic point of the fibre over $\fr{p}$ to $z$, yielding a corresponding chain $C_{\delta}\subset X_{\delta}$. As $\fr{p}$ lies over $\fr{p}_{\delta}$ for each $\delta\in\Delta$, this $C_{\delta}$ is contained in the fibre over $\fr{p}_{\delta}$. 
	\par Therefore, the fibre over $\fr{p}_{\delta}$ has codimension $\geqslant d - h_{\delta}$, as required. This completes the proof.
\es
\bs[Proof of Corollary~\ref{conj:Gabber}]
	Since $\Pic(A)=0$ for any local ring $A$, the claim where $X$ is local follows from \eqref{vanishing:pic}. Consequently, it suffices to prove \eqref{vanishing:pic} and \eqref{vanishing:br}. To do so, we reduce to show that for each integer $n$ \bud\text{$\Pic(X)[n]\iso\Pic(X\setminus Z)[n]$~~~and~~~$H^2(X,\bb{G}_m)[n]\twoheadrightarrow H^2(X\setminus Z,\bb{G}_m)[n]$}.\eud Fix an integer $n$. We first suppose that $d\geqslant 3$. Then, by assumption and Theorem~\ref{thm:intro}, the long exact sequence \bud \cdots\to H^{q-1}_{Z}(X,\mu_n)\to H^{q-1}(X,\mu_n)\to H^{q-1}(X\setminus Z,\mu_n)\to H^{q}_{Z}(X,\mu_n)\to\cdots\eud shows that there is an isomorphism $H^1(X,\mu_n)\cong H^1(X\setminus Z,\mu_n)$. Let us subsequently consider the Kummer exact sequence of sheaves in the $\fppf$ topology on any scheme \bud1\to\mu_n\to\bb{G}_m\xrightarrow{\times n}\bb{G}_m\to 1.\eud Taking the long exact sequence of cohomology groups associated to the above, we obtain a morphism of short exact sequences \bd\label{diag:picard_case} 1\tir& \bb{G}_m(X)/n\arrow[d]\tir&H^1(X,\mu_n)\arrow[d, sloped, "\sim"]\tir&\Pic(X)[n]\arrow[d]\tir&1\\1\tir&\bb{G}_m(X\setminus Z)/n\tir&H^1(X\setminus Z,\mu_n)\tir&\Pic(X\setminus Z)[n]\tir&1,\ed where the middle vertical arrow is an isomorphism. Since $d\geqslant 3$, our hypothesis ensures that the generic fibre of $Z\subset X$ is of codimension $\geqslant 3$, while any non-generic fibre is of codimension $\geqslant 1$. As a consequence, by \ccite{egaIV_3}{théorème~11.3.8 (c)$\implies$(a)}, we deduce that the local ring of $X$ at any point in $Z$ is of depth $\geqslant 2$. Therefore, since $\bb{G}_m$ is an affine scheme, by \ccite{ces-scholze}{Lemma~7.2.7(2)}, it follows that $\bb{G}_m(X)\cong\bb{G}_m(X\setminus Z)$. In particular, we get that the left vertical morphism of \eqref{diag:picard_case} is an isomorphism. Consequently, the argument above shows that $\Pic(X)[n]\iso\Pic(X\setminus Z)[n]$, as required. Thus we are done.
	\par We now assume that $d\geqslant 4$. In this case, by Theorem~\ref{thm:intro} and the long exact sequence of supports, there is an isomorphism $H^2(X,\mu_n)\cong H^2(X\setminus Z,\mu_n)$. Following similar arguments as above, we get a morphism of short exact sequences \bd\label{diag:last} 1\tir&\Pic(X)/n\tir\arrow[d]&H^2(X,\mu_n)\arrow[d, sloped, "\sim"]\tir&H^2(X,\bb{G}_m)[n]\arrow[d]\tir&1\\1\tir&\Pic(X\setminus Z)/n\tir&H^2(X\setminus Z,\mu_n)\tir&H^2(X\setminus Z,\bb{G}_m)[n]\tir&1,\ed where the middle vertical arrow is an isomorphism. This proves that $H^2(X,\bb{G}_m)_{\tors}\twoheadrightarrow H^2(X\setminus Z,\bb{G}_m)_{\tors}$, as required. This concludes the proof.
\es
\br
	In diagram \eqref{diag:last}, we know that $\Pic(A)/n=0$. Therefore, in light of Corollary~\ref{conj:Gabber}, the canonical morphism $\Br(A)[n]\to\Br(U_A)[n]$ is injective if and only if $\Pic(U_A)/n=0$. We expect this to hold when $d\geqslant 4$. In fact, one may hope that $\Pic(U_A)=0$, by extending the arguments of \ccite{sga2}{exposé~XI, \textsection3} to the present setting.
\er
\br[Smooth case]\label{rem:smooth_case}
The arguments in \cite{gabber-ramero-foundations, ning_liu_quasi-split} use different techniques and do not recover the vanishing range established here. In particular, Theorem~\ref{thm:fd_case} appears to be new even when $X$ is smooth over $\mathcal{O}$.
\er
\bs[Proof of Corollary~\ref{conj:bhatt_lurie}]
	The claim \eqref{claim:1:BL} readily follows from Theorem~\ref{thm:intro}. Therefore, it suffices to verify the claim \eqref{claim:2:BL}, for which, we suppose that $d\geqslant 3$.
	\par By \ccite{bhatt_lurie_absolute_prismatic_cohomology}{Proposition~8.4.14 and Variant~8.4.15} and \ccite{brauerIII}{théorème~11.7}, for any scheme $S$, the \textit{syntomic first Chern class map} induces an isomorphism $R\Gamma(S, \mu_{p^r})\iso R\Gamma_{\syn}(S, \bb{Z}/p^r(1))$. In particular, the isomorphism \eqref{claim:2:BL} follows from \eqref{claim:1:BL}. 
\es
\subsection*{Acknowledgements}
The author would like to thank Tess Bouis, Hanlin Cai, Hiroki Kato, Dorin Popescu, and Vaidehee Thatte for many helpful conversations. This project originated during the author’s doctoral studies, and the author is especially grateful to Elden Elmanto and Kęstutis Česnavičius for their many ideas and techniques that significantly influenced this work.
\par This project received partial funding from PNRR grant CF 44/14.11.2022 ``\textit{Cohomological
	Hall algebras of smooth surfaces and applications}'' and NSERC Discovery grant RGPIN-2025-07114, “\textit{Motivic cohomology: theory and applications}”.

\printbibliography{\let\thefootnote\relax\footnote{Simion Stoilow Institute of Mathematics of the Romanian Academy (IMAR), Bucharest, Romania\newline\hspace*{0.48cm} Email: akundu.math@gmail.com.}}{\let\thefootnote\relax\footnote{May 3, 2026}}
\end{document}

%% file: preamble.tex
\usepackage[tmargin=2cm,bmargin=2.5cm,lmargin=1.15cm,rmargin=1.15cm]{geometry}
\usepackage[utf8]{inputenc}
\usepackage[english]{babel}
\usepackage[rm]{titlesec}
\titlelabel{\thetitle.\quad}
\usepackage[T1]{fontenc} 
\usepackage{amssymb,amsfonts}
\usepackage{amsthm}
\usepackage{tikz-cd}
\usepackage[colorinlistoftodos]{todonotes}
\usepackage[framemethod=TikZ]{mdframed} 
\usepackage{calligra,mathrsfs}
\usepackage{amsmath}
\usepackage{mathtools}
\usetikzlibrary{shapes.geometric, arrows}
\usepackage[backend=biber,style=alphabetic,sorting=nyt, maxnames=10, maxalphanames=10]{biblatex}
\usepackage{csquotes}
\usepackage{colonequals}
\usepackage{adjustbox} 
\usepackage[colorlinks, final]{hyperref} 
\usepackage{cleveref}  
\usepackage{enumerate} 
\definecolor{imperialred}{RGB}{237, 41, 57}
\definecolor{royalblue}{RGB}{64, 106, 212}
\definecolor{link}{RGB}{11,0,128}
\hypersetup{
	colorlinks = true,
	breaklinks = true,
	citecolor=royalblue,
	linkcolor=imperialred,
	urlcolor=link,
	linktocpage = true
}
\makeatletter

\makeatother

\newrobustcmd*{\citeauthorfullname}{\AtNextCite{\DeclareNameAlias{labelname}{given-family}}\citeauthor}

\newtheorem{thm}{Theorem}[section]
\newtheorem{cor}[thm]{Corollary}
\newtheorem{prop}[thm]{Proposition}
\newtheorem{lem}[thm]{Lemma}
\newtheorem{conj}[thm]{Conjecture}

\newtheorem*{quest*}{Question}

\theoremstyle{definition}
\newtheorem{defn}[thm]{Definition}

\newtheorem*{exmp*}{Example}

\theoremstyle{remark}

\newtheorem{rem}[thm]{Remark}

\newtheoremstyle{montobbo}
{7pt}
{5pt}
{}
{}%
{\bfseries}
{}%
{.5em}
{\thmnumber{\@{#1}{}~\@{#2}.}%
	\thmnote{~{\bfseries#3.}}}

\theoremstyle{montobbo}
	\newtheorem{montobo}[thm]{}
	\newcommand{\bmo}{\begin{montobo}}
	\newcommand{\emo}{\end{montobo}}
	\newtheorem{talk}{Lecture}
	\newcommand{\btalk}{\begin{talk}}
		\newcommand{\etalk}{\end{talk}}
\usepackage{footnote}
\newcommand{\ccite}[2]{\cite[#2]{#1}}
\newcommand{\stacks}[1]{\ccite{stacks-project}{\href{https://stacks.math.columbia.edu/tag/#1}{Tag \MakeUppercase{#1}}}}	

\newenvironment{manualtheorem}[1]{
	\renewcommand\thethm{#1}
	\thm
}{\endthm}
\newcommand{\but}[1]{\begin{manualtheorem}{#1}}
\newcommand{\eut}{\end{manualtheorem}}
\newenvironment{manualproposition}[1]{
	\renewcommand\thethm{#1}
	\prop
}{\endprop}
\newcommand{\bup}[1]{\begin{manualproposition}{#1}}
	\newcommand{\eup}{\end{manualproposition}}
\newenvironment{manualcorollary}[1]{
	\renewcommand\thethm{#1}
	\cor
}{\endcor}
\newcommand{\buc}[1]{\begin{manualcorollary}{#1}}
	\newcommand{\euc}{\end{manualcorollary}}
	\newenvironment{manualremark}[1]{
		\renewcommand\thethm{#1}
		\rem
	}{\endrem}
	\newcommand{\bur}[1]{\begin{manualremark}{#1}}
		\newcommand{\eur}{\end{manualremark}}
\numberwithin{equation}{thm}
\newenvironment{teqn}{\begin{equation}\begin{tikzcd}\displaystyle}{\end{tikzcd}\end{equation}}
\newenvironment{teqn*}{\begin{equation*}\begin{tikzcd}\displaystyle}{\end{tikzcd}\end{equation*}}
\newenvironment{tpic}{\begin{equation}\begin{tikzpicture}}{\end{tikzpicture}\end{equation}}
\newenvironment{tpic*}{\begin{equation*}\begin{tikzpicture}}{\end{tikzpicture}\end{equation*}}

\DeclareMathAlphabet\matheu{U}{eus}{m}{n}
\DeclareMathAlphabet{\mathmm}{U}{mmambb}{m}{n}
\newcommand{\bb}[1]{\mathbb{#1}}
\newcommand{\fr}[1]{\mathfrak{#1}}
\newcommand{\ca}[1]{\mathcal{#1}}
\newcommand{\scr}[1]{\mathscr{#1}}

\newcommand{\comment}[1]{}
\newcommand{\bun}{\begin{itemize}}
\newcommand{\bn}{\begin{enumerate}}
\newcommand{\bt}{\begin{thm}}
\newcommand{\bl}{\begin{lem}}
\newcommand{\bp}{\begin{prop}}
\newcommand{\bc}{\begin{cor}}
\newcommand{\bd}{\begin{teqn}}
\newcommand{\bud}{\begin{teqn*}}	
\newcommand{\bs}{\begin{proof}}
\newcommand{\br}{\begin{rem}}
\newcommand{\bdf}{\begin{defn}}
\newcommand{\bcj}{\begin{conj}}

\newcommand{\eun}{\end{itemize}}
\newcommand{\en}{\end{enumerate}}
\newcommand{\et}{\end{thm}}
\newcommand{\el}{\end{lem}}
\newcommand{\ep}{\end{prop}}
\newcommand{\ec}{\end{cor}}
\newcommand{\ed}{\end{teqn}}
\newcommand{\eud}{\end{teqn*}}
\newcommand{\es}{\end{proof}}
\newcommand{\er}{\end{rem}}
\newcommand{\edf}{\end{defn}}
\newcommand{\ecj}{\end{conj}}





\DeclareMathOperator{\colim}{colim}

\DeclareMathOperator{\spec}{Spec}

\DeclareMathOperator{\Pic}{Pic}



\DeclareMathOperator{\fppf}{\textit{fppf}}

\DeclareMathOperator{\syn}{syn}
\DeclareMathOperator{\tors}{tors}
\DeclareMathOperator{\Br}{Br}


\newcommand{\tir}{\arrow[r]}

\newcommand{\iso}{\xrightarrow{\sim}}


\makeatletter\newcommand{\diso}{\arrow[r, "\sim"]}\makeatother

%% file: bibli.bib
@article {antieau_datta_derived_splinters,
	AUTHOR = {Antieau, Benjamin and Datta, Rankeya},
	TITLE = {Valuation rings are derived splinters},
	JOURNAL = {Math. Z.},
	FJOURNAL = {Mathematische Zeitschrift},
	VOLUME = {299},
	YEAR = {2021},
	NUMBER = {1-2},
	PAGES = {827--851},
	ISSN = {0025-5874,1432-1823},
	MRCLASS = {14F08 (13A18 13D09 13D22 14B05)},
	MRNUMBER = {4311620},
	MRREVIEWER = {Igor\ A.\ Rapinchuk},
	DOI = {10.1007/s00209-020-02683-6},
	URL = {https://doi.org/10.1007/s00209-020-02683-6},
}

@misc{bhatt_lurie_absolute_prismatic_cohomology,
	title={Absolute prismatic cohomology}, 
	author={Bhargav Bhatt and Jacob Lurie},
	year={2022},
	eprint={2201.06120},
	archivePrefix={arXiv},
	primaryClass={math.AG},
	url={https://arxiv.org/abs/2201.06120}, 
}

@article {bhatt-mathew-arc,
    AUTHOR = {Bhatt, Bhargav and Mathew, Akhil},
     TITLE = {The arc-topology},
   JOURNAL = {Duke Math. J.},
  FJOURNAL = {Duke Mathematical Journal},
    VOLUME = {170},
      YEAR = {2021},
    NUMBER = {9},
     PAGES = {1899--1988},
      ISSN = {0012-7094},
   MRCLASS = {14F20 (14F06 14G22)},
  MRNUMBER = {4278670},
       DOI = {10.1215/00127094-2020-0088},
}

@article{ces-brauer,
    AUTHOR = {Česnavičius, Kęstutis},
     TITLE = {Purity for the {B}rauer group},
   JOURNAL = {Duke Math. J.},
  FJOURNAL = {Duke Mathematical Journal},
    VOLUME = {168},
      YEAR = {2019},
    NUMBER = {8},
     PAGES = {1461--1486},
      ISSN = {0012-7094},
   MRCLASS = {14F22 (14F20 14G22 16K50)},
  MRNUMBER = {3959863},
MRREVIEWER = {J\"{o}rg Jahnel},
       DOI = {10.1215/00127094-2018-0057},
}

@article {ces-scholze,
    AUTHOR = {Česnavičius, Kęstutis and Scholze, Peter},
     TITLE = {Purity for flat cohomology},
   JOURNAL = {Ann. of Math. (2)},
  FJOURNAL = {Annals of Mathematics. Second Series},
    VOLUME = {199},
      YEAR = {2024},
    NUMBER = {1},
     PAGES = {51--180},
      ISSN = {0003-486X,1939-8980},
   MRCLASS = {14F20 (14F22 14F30 14H20 18G90)},
  MRNUMBER = {4681144},
       DOI = {10.4007/annals.2024.199.1.2},
       URL = {https://doi.org/10.4007/annals.2024.199.1.2},
}

@book {ct-skorobogatov_brauer_grothendieck_group_book,
    AUTHOR = {Colliot-Th\'{e}l\`ene, Jean-Louis and Skorobogatov, Alexei N.},
     TITLE = {The {B}rauer-{G}rothendieck group},
    SERIES = {Ergebnisse der Mathematik und ihrer Grenzgebiete. 3. Folge. A
              Series of Modern Surveys in Mathematics [Results in
              Mathematics and Related Areas. 3rd Series. A Series of Modern
              Surveys in Mathematics]},
    VOLUME = {71},
 PUBLISHER = {Springer, Cham},
      YEAR = {[2021] \copyright 2021},
     PAGES = {xv+453},
      ISBN = {978-3-030-74247-8; 978-3-030-74248-5},
   MRCLASS = {14F22 (14E08 14G05 14G12 14K05)},
  MRNUMBER = {4304038},
MRREVIEWER = {Thomas Benedict Williams},
       DOI = {10.1007/978-3-030-74248-5},
}

@article {datta_local_cohomology_valuation_ring,
	AUTHOR = {Datta, Rankeya},
	TITLE = {({N}on)vanishing results on local cohomology of valuation
	rings},
	JOURNAL = {J. Algebra},
	FJOURNAL = {Journal of Algebra},
	VOLUME = {479},
	YEAR = {2017},
	PAGES = {413--436},
	ISSN = {0021-8693,1090-266X},
	MRCLASS = {13F30 (13D45)},
	MRNUMBER = {3627291},
	MRREVIEWER = {Franz-Viktor\ Kuhlmann},
	DOI = {10.1016/j.jalgebra.2016.12.032},
	URL = {https://doi.org/10.1016/j.jalgebra.2016.12.032},
}

@misc{de_jong_gabber_br,
 	title       = {A result of Gabber},
 	author      = {de Jong, A. J.},
 	url         = {https://www.math.columbia.edu/~dejong/papers/2-gabber.pdf},
 	shorthand= {dJ02},
 }

@article {egaIV_3,
    shorthand={EGA IV$_3$},
    AUTHOR = {Grothendieck, A.},
     TITLE = {\'{E}l\'{e}ments de g\'{e}om\'{e}trie alg\'{e}brique. {IV}. \'{E}tude locale des
              sch\'{e}mas et des morphismes de sch\'{e}mas. {III}},
   JOURNAL = {Inst. Hautes \'{E}tudes Sci. Publ. Math.},
  FJOURNAL = {Institut des Hautes \'{E}tudes Scientifiques. Publications
              Math\'{e}matiques},
    NUMBER = {28},
      YEAR = {1966},
     PAGES = {255},
      ISSN = {0073-8301},
   MRCLASS = {14.55},
  MRNUMBER = {217086},
MRREVIEWER = {J. P. Murre},
       URL = {http://www.numdam.org/item?id=PMIHES_1966__28__255_0},
}

@article {ehik,
    AUTHOR = {Elmanto, Elden and Hoyois, Marc and Iwasa, Ryomei and Kelly,
              Shane},
     TITLE = {Cdh descent, cdarc descent, and {M}ilnor excision},
   JOURNAL = {Math. Ann.},
  FJOURNAL = {Mathematische Annalen},
    VOLUME = {379},
      YEAR = {2021},
    NUMBER = {3-4},
     PAGES = {1011--1045},
      ISSN = {0025-5831},
   MRCLASS = {14F42 (14A20)},
  MRNUMBER = {4238259},
       DOI = {10.1007/s00208-020-02083-5},
}

@book {fujiwara_kato_foundations_rigid_geometry,
    AUTHOR = {Fujiwara, Kazuhiro and Kato, Fumiharu},
     TITLE = {Foundations of rigid geometry. {I}},
    SERIES = {EMS Monographs in Mathematics},
 PUBLISHER = {European Mathematical Society (EMS), Z\"{u}rich},
      YEAR = {2018},
     PAGES = {xxxiv+829},
      ISBN = {978-3-03719-135-4},
   MRCLASS = {14G22 (14F05)},
  MRNUMBER = {3752648},
MRREVIEWER = {Christopher David Lazda},
}

@incollection {gabber_oberwolfach_report,
	TITLE = {Arithmetic algebraic geometry},
	NOTE = {Abstracts from the workshop held August 1--7, 2004,
	Organized by Gerd Faltings, G\"unter Harder and Nicholas M.
	Katz,
	Oberwolfach Reports. Vol. 1, no. 3},
	JOURNAL = {Oberwolfach Rep.},
	FJOURNAL = {Oberwolfach Reports},
	VOLUME = {1},
	YEAR = {2004},
	NUMBER = {3},
	PAGES = {1971--2013},
	ISSN = {1660-8933,1660-8941},
	MRCLASS = {11-06 (14-06 22-06)},
	MRNUMBER = {2144155},
	DOI = {10.4171/OWR/2004/37},
	URL = {https://doi.org/10.4171/OWR/2004/37},
	shorthand= {Gab04},
}

@article {gabber_kelly_points,
	AUTHOR = {Gabber, Ofer and Kelly, Shane},
	TITLE = {Points in algebraic geometry},
	JOURNAL = {J. Pure Appl. Algebra},
	FJOURNAL = {Journal of Pure and Applied Algebra},
	VOLUME = {219},
	YEAR = {2015},
	NUMBER = {10},
	PAGES = {4667--4680},
	ISSN = {0022-4049,1873-1376},
	MRCLASS = {14F05},
	MRNUMBER = {3346512},
	MRREVIEWER = {Maria\ Chiara\ Brambilla},
	DOI = {10.1016/j.jpaa.2015.03.001},
	URL = {https://doi.org/10.1016/j.jpaa.2015.03.001},
}

@book {gabber-ramero-almost,
    AUTHOR = {Gabber, Ofer and Ramero, Lorenzo},
     TITLE = {Almost ring theory},
    SERIES = {Lecture Notes in Mathematics},
    VOLUME = {1800},
 PUBLISHER = {Springer-Verlag, Berlin},
      YEAR = {2003},
     PAGES = {vi+307},
      ISBN = {3-540-40594-1},
   MRCLASS = {13D10 (13B40 13D03 14G22 18D10)},
  MRNUMBER = {2004652},
       DOI = {10.1007/b10047},
}

@misc{gabber-ramero-foundations,
    title={Foundations for almost ring theory -- Release 7.5},
    author={Ofer Gabber and Lorenzo Ramero},
    shorthand={GR18},
    year={2004},
    eprint={math/0409584},
    archivePrefix={arXiv},
    primaryClass={math.AG}
}

@incollection {brauerIII,
    AUTHOR = {Grothendieck, Alexander},
     TITLE = {Le groupe de {B}rauer. {III}. {E}xemples et compl\'{e}ments},
 BOOKTITLE = {Dix expos\'{e}s sur la cohomologie des sch\'{e}mas},
    SERIES = {Adv. Stud. Pure Math.},
    VOLUME = {3},
     PAGES = {88--188},
 PUBLISHER = {North-Holland, Amsterdam},
      YEAR = {1968},
   MRCLASS = {14F22},
  MRNUMBER = {244271},
MRREVIEWER = {James Milne},
}

@book {huber96,
    AUTHOR = {Huber, Roland},
     TITLE = {\'{E}tale cohomology of rigid analytic varieties and adic spaces},
    SERIES = {Aspects of Mathematics, E30},
 PUBLISHER = {Friedr. Vieweg \& Sohn, Braunschweig},
      YEAR = {1996},
     PAGES = {x+450},
      ISBN = {3-528-06794-2},
   MRCLASS = {14G22 (14F20)},
  MRNUMBER = {1734903},
MRREVIEWER = {Lorenzo Ramero},
       DOI = {10.1007/978-3-663-09991-8},
}

@article {kato_thatte_upper_Ramification_valuation,
	AUTHOR = {Kato, Kazuya and Thatte, Vaidehee},
	TITLE = {Upper ramification groups for arbitrary valuation rings},
	JOURNAL = {Tunis. J. Math.},
	FJOURNAL = {Tunisian Journal of Mathematics},
	VOLUME = {6},
	YEAR = {2024},
	NUMBER = {4},
	PAGES = {589--646},
	ISSN = {2576-7658,2576-7666},
	MRCLASS = {11G99 (14G22)},
	MRNUMBER = {4843477},
	DOI = {10.2140/tunis.2024.6.589},
	URL = {https://doi.org/10.2140/tunis.2024.6.589},
}

@article {kelly-morrow,
    AUTHOR = {Kelly, Shane and Morrow, Matthew},
     TITLE = {{$K$}-theory of valuation rings},
   JOURNAL = {Compos. Math.},
  FJOURNAL = {Compositio Mathematica},
    VOLUME = {157},
      YEAR = {2021},
    NUMBER = {6},
     PAGES = {1121--1142},
      ISSN = {0010-437X},
   MRCLASS = {19D50 (14F30)},
  MRNUMBER = {4264079},
       DOI = {10.1112/s0010437x21007119},
}

@article {kerz-strunk-tamme,
    AUTHOR = {Kerz, Moritz and Strunk, Florian and Tamme, Georg},
     TITLE = {Towards {V}orst's conjecture in positive characteristic},
   JOURNAL = {Compos. Math.},
  FJOURNAL = {Compositio Mathematica},
    VOLUME = {157},
      YEAR = {2021},
    NUMBER = {6},
     PAGES = {1143--1171},
      ISSN = {0010-437X},
   MRCLASS = {19D35 (14F42 19D55)},
  MRNUMBER = {4270122},
       DOI = {10.1112/S0010437X21007120},
}

@misc{amar:unramified_case_gersten_conjecture,
      title={Gersten's Injectivity for Smooth Algebras over Valuation Rings}, 
      author={Arnab Kundu},
      year={2024},
      eprint={2404.06655},
      archivePrefix={arXiv},
      primaryClass={math.AG},
      url={https://arxiv.org/abs/2404.06655},
}

@unpublished{phd-thesis,
    title={Torsors on Smooth Algebras over Valuation Rings},
    author={Kundu, Arnab},
    year={2023},
    type={PhD Thesis},
    url={https://cnrs.hal.science/tel-04129990/},
    keywords={amar},
    %shorthand={Thesis},
}

@misc{madapusi_mondal_perfect_f_gauges,
	title={Perfect $F$-gauges and finite flat group schemes}, 
	author={Keerthi Madapusi and Shubhodip Mondal},
	year={2025},
	eprint={2509.01573},
	archivePrefix={arXiv},
	primaryClass={math.NT},
	url={https://arxiv.org/abs/2509.01573}, 
}

@article{ning_liu_quasi-split,
    AUTHOR = {Guo, Ning and Liu, Fei},
     TITLE = {Purity and quasi-split torsors over {P}r\"ufer bases},
   JOURNAL = {J. \'Ec. polytech. Math.},
  FJOURNAL = {Journal de l'\'Ecole polytechnique. Math\'ematiques},
    VOLUME = {11},
      YEAR = {2024},
     PAGES = {187--246},
      ISSN = {2429-7100,2270-518X},
   MRCLASS = {14F20 (13F05 14F22 14G22 16K50)},
  MRNUMBER = {4695964},
       DOI = {10.5802/jep.253},
       URL = {https://doi.org/10.5802/jep.253},
}

@article{popescu_neron,
    AUTHOR = {Popescu, Dorin},
    TITLE = {N\'eron desingularization of extensions of valuation rings},
    BOOKTITLE = {Transcendence in algebra, combinatorics, geometry and number
    theory},
    SERIES = {Springer Proc. Math. Stat.},
    VOLUME = {373},
    PAGES = {275--307},
    NOTE = {With an appendix by Kęstutis Česnavičius},
    PUBLISHER = {Springer, Cham},
    YEAR = {[2021] \copyright 2021},
    ISBN = {978-3-030-84303-8; 978-3-030-84304-5},
    MRCLASS = {13F30 (03C20 13B02 14B05)},
    MRNUMBER = {4378980},
    MRREVIEWER = {Lokendra\ Paudel},
    }

@article {popescu_algebraic_complete_intersection,
    	AUTHOR = {Popescu, Dorin},
    	TITLE = {Algebraic valuation ring extensions as limits of complete
    	intersection algebras},
    	JOURNAL = {Rev. Mat. Complut.},
    	FJOURNAL = {Revista Matem\'atica Complutense},
    	VOLUME = {37},
    	YEAR = {2024},
    	NUMBER = {2},
    	PAGES = {467--472},
    	ISSN = {1139-1138,1988-2807},
    	MRCLASS = {13F30 (13A18 13B40 13F20)},
    	MRNUMBER = {4739543},
    	MRREVIEWER = {Michael\ Steward},
    	DOI = {10.1007/s13163-023-00468-z},
    	URL = {https://doi.org/10.1007/s13163-023-00468-z},
    }

@article {popescu_free_value_group_complete_intersection,
    	AUTHOR = {Popescu, Dorin},
    	TITLE = {Filtered colimits of complete intersection algebras},
    	JOURNAL = {Rev. Roumaine Math. Pures Appl.},
    	FJOURNAL = {Revue Roumaine de Math\'ematiques Pures et Appliqu\'ees.
    	Romanian Journal of Pure and Applied Mathematics},
    	VOLUME = {70},
    	YEAR = {2025},
    	NUMBER = {1-2},
    	PAGES = {157--165},
    	ISSN = {0035-3965},
    	MRCLASS = {13F30 (13A18 13B40 13F20)},
    	MRNUMBER = {4876699},
    	MRREVIEWER = {Carmelo\ Antonio\ Finocchiaro},
    }

@misc{popescu_char_p_complete_intersection,
  		title={Valuation rings as limits of complete intersection rings}, 
  		author={Dorin Popescu},
  		year={2026},
  		eprint={2004.11004},
  		archivePrefix={arXiv},
  		primaryClass={math.AC},
  		url={https://arxiv.org/abs/2004.11004}, 
  	}

@misc{popescu_mixed_char_torsion_free_value_group_complete_intersection,
    	title={Valuation rings of mixed characteristic as limits of complete intersection rings}, 
    	author={Dorin Popescu},
    	year={2026},
    	eprint={2101.01346},
    	archivePrefix={arXiv},
    	primaryClass={math.AC},
    	url={https://arxiv.org/abs/2101.01346}, 
    }

@article {popescu_dimension_1_free_value_group_smooth,
    	AUTHOR = {Popescu, Dorin},
    	TITLE = {Valuation rings of dimension one as limits of smooth algebras},
    	JOURNAL = {Bull. Math. Soc. Sci. Math. Roumanie (N.S.)},
    	FJOURNAL = {Bulletin Math\'ematique de la Soci\'et\'e{} des Sciences
    	Math\'ematiques de Roumanie. Nouvelle S\'erie},
    	VOLUME = {64(112)},
    	YEAR = {2021},
    	NUMBER = {1},
    	PAGES = {63--73},
    	ISSN = {1220-3874,2065-0264},
    	MRCLASS = {13F30 (13A18 13B40 13L05)},
    	MRNUMBER = {4243507},
    	MRREVIEWER = {Chudamani\ Pranesachar\ Anil Kumar},
    }

@article {popescu_immediate_pure_transcendental_smooth,
    	AUTHOR = {Popescu, Dorin},
    	TITLE = {Pure transcendental, immediate valuation ring extensions as
    	limits of smooth algebras},
    	JOURNAL = {Manuscripta Math.},
    	FJOURNAL = {Manuscripta Mathematica},
    	VOLUME = {176},
    	YEAR = {2025},
    	NUMBER = {5},
    	PAGES = {Paper No. 63, 22},
    	ISSN = {0025-2611,1432-1785},
    	MRCLASS = {13F30 (13A18 13B40 13F20)},
    	MRNUMBER = {4949577},
    	DOI = {10.1007/s00229-025-01635-w},
    	URL = {https://doi.org/10.1007/s00229-025-01635-w},
    }

@book {sga2,
	AUTHOR = {Grothendieck, Alexander},
	TITLE = {Cohomologie locale des faisceaux coh\'erents et th\'eor\`emes
	de {L}efschetz locaux et globaux ({SGA} 2)},
	shorthand = {SGA~$2$},
	SERIES = {Documents Math\'ematiques (Paris) [Mathematical Documents
	(Paris)]},
	VOLUME = {4},
	EDITOR = {Laszlo, Yves},
	NOTE = {S\'eminaire de G\'eom\'etrie Alg\'ebrique du Bois Marie, 1962.,
	Augment\'e{} d'un expos\'e{} de Mich\`ele Raynaud. [With an
	expos\'e{} by Mich\`ele Raynaud],
	Revised reprint of the 1968 French original},
	PUBLISHER = {Soci\'et\'e{} Math\'ematique de France, Paris},
	YEAR = {2005},
	PAGES = {x+208},
	ISBN = {2-85629-169-4},
	MRCLASS = {14B15 (14C20 14F20)},
	MRNUMBER = {2171939},
}

@book {sga4ii,
     TITLE = {Th\'{e}orie des topos et cohomologie \'{e}tale des sch\'{e}mas. {T}ome 2},
     shorthand = {SGA~$4_{\text{II}}$},
    SERIES = {Lecture Notes in Mathematics, Vol. 270},
      NOTE = {S\'{e}minaire de G\'{e}om\'{e}trie Alg\'{e}brique du Bois-Marie 1963--1964
              (SGA 4),
              Dirig\'{e} par M. Artin, A. Grothendieck et J. L. Verdier. Avec la
              collaboration de N. Bourbaki, P. Deligne et B. Saint-Donat},
 PUBLISHER = {Springer-Verlag, Berlin-New York},
      YEAR = {1972},
     PAGES = {iv+418},
   MRCLASS = {14-06},
  MRNUMBER = {0354653},
}

@misc{stacks-project,
	author       	= {The {Stacks project authors}},
	title        	= {The Stacks project},
	howpublished 	= {\url{https://stacks.math.columbia.edu}},
	year         	= {2026},
}

@article {zariski_local_uniformisation,
	AUTHOR = {Zariski, Oscar},
	TITLE = {Local uniformization on algebraic varieties},
	JOURNAL = {Ann. of Math. (2)},
	FJOURNAL = {Annals of Mathematics. Second Series},
	VOLUME = {41},
	YEAR = {1940},
	PAGES = {852--896},
	ISSN = {0003-486X},
	MRCLASS = {09.1X},
	MRNUMBER = {2864},
	MRREVIEWER = {R.\ J.\ Walker},
	DOI = {10.2307/1968864},
	URL = {https://doi.org/10.2307/1968864},
}
